\documentclass{article}

\usepackage{amssymb}
\usepackage{amsmath}
\usepackage{tensor}
\newcommand{\proof}{{\bf Proof:  }}
\newcommand{\remark}{{\bf Remark:  }}

\newcommand{\example}{{\bf Example:  }}

\newcommand{\hb}{\newline\hspace*{\fill}$\Box$}

\newtheorem{theorem}{Theorem}[section]
\newtheorem{lemma}[theorem]{Lemma}
\newtheorem{definition}[theorem]{Definition}
\newtheorem{proposition}[theorem]{Proposition}
\newtheorem{corollary}[theorem]{Corollary}

\newtheorem{conjecture}[theorem]{Conjecture}

\title{Degenerate Cohomological Hall algebra and quantized Donaldson-Thomas invariants for $m$-loop quivers}
\author{Markus Reineke}

\begin{document}

\parindent0pt

\maketitle

\begin{abstract} We derive a combinatorial formula for quantized Donaldson-Thomas invariants of the $m$-loop quiver. Our main tools are the combinatorics of noncommutative Hilbert schemes and a degenerate version of the Cohomological Hall algebra of this quiver.
\end{abstract}

\noindent

\section{Introduction}

Generalized Donaldson-Thomas invariants of (noncommutative) varieties arise from factorizations of generating series of motivic invariants of Hilbert scheme-type varieties into Euler products. For $3$-Calabi-Yau manifolds, this principle is developed extensively in \cite{KSDT1}.\\[1ex]
In \cite{RePois}, the author showed that the wall-crossing formulae of \cite{KSDT1} can be modeled using Hilbert schemes of path algebras of quivers; explicit calculations for these varieties in \cite{ReFE} allowed to derive relative integrality (that is, preservation of integrality under wall-crossing) of generalized Donaldson-Thomas invariants.\\[1ex]
In \cite{KSCOHA} a general framework for the study of such integrality properties is proposed, the central tools being Cohomological Hall algebras and the geometric concept of factorization systems.\\[2ex]
The purpose of the present paper is to develop an explicit, in most parts purely combinatorial, setup for the study of the quantized Donaldson-Thomas invariants of \cite{KSCOHA} in the very special, but typical, case of the $m$-loop quiver. The relevant concepts of \cite{KSCOHA} are discussed in sections \ref{dt1}, \ref{coha}. Our approach is based on the explicit description of Hilbert schemes attached to this quiver of \cite{ReHilb1}, which is reviewed in Section \ref{Hilb}. It allows us to give a combinatorial description of a degenerate version of the Cohomological Hall algebra, whose structure is easily described (see Section \ref{dcoha}). Using number-theoretic arguments similar to \cite{ReFE}, we obtain explicit formulas for these quantized Donaldson-Thomas invariants (see Theorem \ref{mainresult}) in terms of cyclic classes of certain integer sequences in Section \ref{integrality}. We also relate this combinatorics to a similar one appearing in the study of Higgs moduli in \cite{HT}, see Section \ref{higgs}.\\[1ex]
Roughly spoken, our approach uses (the combinatorics of) noncommutative Hilbert schemes as a transitional tool between the geometric problem of determination of Donaldson-Thomas invariants and the combinatorial object of cyclic configurations. However, the present approach is not strong enough to yield the positivity properties conjectured in \cite{KSCOHA}.\\[3ex]
{\bf Acknowledgments:} The author would like to thank S. Mozgovoy for several discussions concerning this work, and in particular for pointing out the potential relation to Higgs moduli. This work was started while the author participated in the workshop ``Wall-crossing in Mathematics and Physics'' at Urbana-Champaign, and finished during stays at the Issac Newton Institute Cambridge and the Hausdorff Institute Bonn. The author would like to thank the organizers and participants of these programmes for the inspiring atmosphere. 

\section{Noncommutative Hilbert schemes}\label{Hilb}

In this section, we recall the definition of noncommutative Hilbert schemes and their main properties following \cite{ReHilb1}. We also relate the relevant combinatorics of trees to a combinatorics of partitions which will play a major role in the following.\\[2ex]
Fix an integer $m\geq 1$. For $n\geq 0$, we call a pair consisting of a tuple $(\varphi_1,\ldots,\varphi_m)$ of linear operators on ${\bf C}^n$ and a vector $v\in {\bf C}^n$ stable if $v$ is cyclic for the representation of the free algebra $F^{(m)}={\bf C}\langle x_1,\ldots,x_m\rangle$ on ${\bf C}^n$ defined by the operators $\varphi_i$, that is, if ${\bf C}\langle \varphi_1,\ldots,\varphi_m\rangle v={\bf C}^n$. This defines an open subset of the affine space ${\rm End}({\bf C}^n)^m\oplus {\bf C}^n$, for which a geometric quotient by the action of ${\rm GL}_n({\bf C})$ via $g(\varphi_1,\ldots,\varphi_m,v)=(g\varphi_1g^{-1},\ldots,g\varphi_mg^{-1},gv)$ exists. This quotient is denoted by ${\rm Hilb}^{(m)}_n$ and is called a noncommutative Hilbert scheme for $F^{(m)}$: in analogy with the Hilbert scheme of $n$ points of an affine variety $X$ parametrizing codimension $n$ ideals in the coordinate ring of $X$, the variety ${\rm Hilb}_n^{(m)}$ parametrizes left ideals $I$ in $F^{(m)}$ of codimension $n$, that is, ideals such that $\dim_{\bf C}F^{(m)}/I=n$. Namely, to a tuple as above we associate the left ideal of polynomials $P\in F^{(m)}$ such that $F(\varphi_1,\ldots,\varphi_m)v=0$. Conversely, given a left ideal $I\subset F^{(m)}$, we choose an isomorphism between $F^{(m)}/I$ and ${\bf C}^n$. The operators $\varphi_i$ are induced by left multiplication on $F^{(m)}/I$ via this isomorphism, whereas $v$ is induced by the coset of the unit $1\in F^{(m)}$. This tuple is stable by definition, and well defined up to the choice of the isomorphism $F^{(m)}/I\simeq{\bf C}^n$, that is, up to the ${\rm GL}_n({\bf C})$-action.\\[1ex]
Consider the set $\Omega^{(m)}$ of words $\omega=i_1\ldots i_k$ in the alphabet $\{1,...m\}$. Composition of words defines a monoid structure on $\Omega^{(m)}$; define $\omega'$ to be a left subword of $\omega$ if $\omega=\omega'\omega''$ for a word $\omega''$. The set $\Omega^{(m)}$ carries a lexicographic ordering $\leq_{\rm lex}$ induced by the canonical total ordering on the alphabet $\{1,\ldots,m\}$. An $m$-ary tree is a subset $T\subset\Omega^{(m)}$ which is closed under left subwords. This terminology is explained as follows: a subset $T$ is visualized as the tree with nodes $\omega$ for $\omega\in T$ and an edge of colour $i$ from $\omega$ to $\omega i$ if $\omega,\omega i\in T$; the empty word corresponds to the root of the tree.\\[1ex]
For a tree $T$, define its corona $C(T)$ as the set of all $\omega\in\Omega^{(m)}$ such that $\omega\not\in T$, but $\omega'\in T$ for $\omega=\omega'i$. We have $|C(T)|=(m-1)|T|+1$.\\[1ex]
Given a word $\omega=i_1\ldots i_k$ and a tuple of operators $(\varphi_1,\ldots,\varphi_m)$ as above, we define $\varphi_\omega=\varphi_{i_k}\circ\ldots\circ\varphi_{i_1}$. For a tree $T$ of cardinality $n$, define $Z_T\subset{\rm Hilb}^{(m)}_n$ as the set of classes of tuples $(\varphi_1,\ldots,\varphi_m,v)$ such that:
\begin{enumerate}
\item the elements $\varphi_\omega v$ for $\omega\in T$ form a basis of $V$,
\item if $\omega\in C(T)$, then $\varphi_\omega v=\sum_{\omega'}\lambda_{\omega,\omega'}v$, where the sum ranges over all words $\omega'\in T$ such that $\omega'<_{\rm lex}\omega$.
\end{enumerate}
Denote by $d(T)$ the number of pairs $(\omega,\omega')$ such that $\omega\in C(T)$, $\omega'\in T$ and $\omega'<_{\rm lex}\omega$.

\begin{theorem}\cite[Theorem 1.3]{ReHilb1} The following holds:
\begin{enumerate}
\item $Z_T$ is a locally closed subset of ${\rm Hilb}_n^{(m)}$, which is isomorphic to an affine space of dimension $d(T)$.
\item The subsets $Z_T$, for $T$ ranging over all trees of cardinality $n$, define a cell decomposition of ${\rm Hilb}_n^{(m)}$, that is, there exists a decreasing filtration of ${\rm Hilb}_n^{(m)}$ by closed subvarieties, such that the successive complements are the subsets $Z_T$.
\end{enumerate}
\end{theorem}

As a corollary to this geometric description, we can derive precise information on the cohomology (singular cohomology with rational coefficients) of ${\rm Hilb}_n^{(m)}$. The existence of a cell decomposition implies vanishing of odd cohomology (and algebraicity of even cohomology), thus we can consider the following generating series of Poincar\'e polynomials
$$F(q,t)=\sum_{n\geq 0}q^{(m-1)\binom{n}{2}}\sum_k\dim H^k({\rm Hilb}_n^{(m)})q^{-k/2}t^n\in{\bf Z}[q,q^{-1}][[t]],$$
as well as its specialization
$$F(t)=F(1,t)=\sum_{n\geq 0}\chi({\rm Hilb}_n^{(m)})t^n\in{\bf Z}[[t]].$$
We also define
$$H(q,t)=\sum_{n\geq 0}\frac{q^{(m-1)\binom{n}{2}}}{(1-q^{-1})\cdot\ldots\cdot(1-q^{-n})}t^n\in{\bf Q}(q)[[t]],$$
which is a $q$-hypergeometric series whose major role for the following will be explained in the next section.

\begin{corollary} We have the following explicit descriptions of the series $F(q,t)$ and $F(t)$:
\begin{enumerate}
\item The series $F(q,t)$ is uniquely determined as the solution in ${\bf Q}(q)[[t]]$ of the algebraic $q$-difference equation
$$F(q,t)=1+t\prod_{k=0}^{m-1}F(q,q^kt).$$
\item The series $F(t)$ is uniquely determined as the solution in ${\bf Q}[[t]]$ of the algebraic equation
$$F(t)=1+tF(t)^m.$$
\item The Euler characteristic of ${\rm Hilb}_n^{(m)}$ equals the number of $m$-ary trees with $n$ nodes, which is $\frac{1}{(m-1)n+1}{\binom{mn}{n}}$.
\item We have $F(q,t)=\frac{H(q,t)}{H(q,q^{-1}t)}$.
\end{enumerate}
\end{corollary}

\proof In the notation of \cite{ReHilb1}, the series $F(q,t)$ equals the series $\overline{\zeta}^{(m)}_1(q,t)$ of \cite[Section 5]{ReHilb1} by \cite[Corollary 4.4.]{ReHilb1}. The first statement translates the operation of grafting of trees; see \cite[Theorem 5.5.]{ReHilb1}. Specialization of the functional equation to $q=1$ yields the second statement. The third statement follows from an explicit formula for the number of $m$-ary trees; see \cite[Corollary 4.5.]{ReHilb1}. The fourth statement is a special case of \cite[Theorem 5.2.]{ERHilb2}; in the present case, it is easily derived from the identity
$$H(q,t)=H(q,q^{-1}t)+tH(q,q^{m-1}t)$$
(which follows by a direct calculation from the definition of $H(q,t)$), together with the first statement.\hb

Denote by $T_n$ the set of partitions $0\leq \lambda_1\leq\ldots\leq \lambda_n$ such that $\lambda_i\leq(m-1)(i-1)$ for all $i=1,\ldots,n$. Define the weight of $\lambda\in T_n$ as ${\rm wt}(\lambda)=(m-1)\binom{n}{2}-|\lambda|$. We also define a weight function ${\rm wt}(T)$ on trees $T$ as above by
$${\rm wt}(T)=(m-1){\binom{|T|}{2}}-|\{(\omega',\omega)\in C(T)\times T\,:\, \omega'<_{\rm lex}\omega\}|,$$
thus ${\rm wt}(T)=d(T)-(m-2){\binom{n+1}{2}}-n$ by definition of $d(T)$.\\[1ex]
Given an $m$-ary tree $T\subset\Omega$ with $n$ vertices as above, write $T=\{\omega_1,\ldots,\omega_n\}$ with $\omega_1<_{\rm lex}\ldots <_{\rm lex}\omega_n$. We define a partition $\lambda(T)$ by $$\lambda(T)_i=|\{\omega\in C(T)\,:\, \omega<_{\rm lex}\omega_i\}|.$$

\begin{proposition} The map associating $\lambda(T)$ to $T$ defines a weight-preserving bijection between $m$-ary trees with $n$ nodes and $T_n$.
\end{proposition}

\proof To prove that $\lambda(T)$ belongs to $T_n$, we note that an element $\omega\in C(T)$ such that $\omega<_{\rm lex}\omega_k$ belongs to $C(T_k)\setminus\{\omega_k\}$ for the subtree $T_k=\{\omega_1,\ldots,\omega_{k-1}\}$ of $T$; this is a set of cardinality $(m-1)(k-1)$. We reconstruct the tree from the partition $\lambda\in T_n$ inductively as follows: we start with the empty tree $T_0$. In the $k$-th step, we list the elements of the corona of $T_{k-1}$ in ascending lexicographic order as $C(T_{k-1})=\{\omega^k_i,\ldots,\omega^k_{(m-1)(k-1)-1}\}$ and define $T_k=T_{k-1}\cup\{\omega^k_{\lambda_k+1}\}$. We then have $\{\omega\in C(T),\:\, \omega<_{\rm lex}\omega_k\}=\{\omega^k_1,\ldots,\omega^k_{\lambda_k}\}$, proving that $T$ is reconstructed from $\lambda(T)$. The equality of the weights of $T$ and $\lambda(T)$ follows from the definitions.\hb

\section{Donaldson-Thomas type invariants}\label{dt1}

The following definition of Donaldson-Thomas type invariants for the $m$-loop quiver is motivated by \cite{KSDT1}.

\begin{definition} Define ${\rm DT}_n^{(m)}\in{\bf Q}$ for $n\geq 1$ by writing $$F((-1)^{m-1}t)=\prod_{n\geq 1}(1-t^n)^{-(-1)^{(m-1)n}n{\rm DT}^{(m)}_n}.$$ 
\end{definition}

These numbers are well-defined since $F(t)$ is an integral power series with constant term $1$. A priori, we have $n{\rm DT}^{(m)}_n\in{\bf Z}$.

\begin{theorem}\label{int}\cite{ReFE} We have ${\rm DT}^{(m)}_n\in{\bf N}$; explicitly, these numbers are given by the following formula:
$${\rm DT}_n^{(m)}=\frac{1}{n^2}\sum_{d|n}\mu(\frac{n}{d})(-1)^{(m-1)(n-d)}{\binom{mn-1}{n-1}}.$$
\end{theorem}

We make this formula more explicit by giving some examples; note that ${\rm DT}_n^{(m)}$ is a polynomial in $m$ except if $n\equiv 2\bmod 4$ (this phenomenon will become more transparent in the following sections).

$${\rm DT}^{(m)}_1=1,\;\;\;
{\rm DT}^{(m)}_2=\left\lfloor \frac{m}{2}\right\rfloor,\;\;\;
{\rm DT}^{(m)}_3=\frac{m(m-1)}{2},$$
$${\rm DT}^{(m)}_4=\frac{m(m-1)(2m-1)}{3},\;\;\;
{\rm DT}^{(m)}_5=\frac{5m(m-1)(5m^2-5m+2)}{24},$$
$${\rm DT}^{(m)}_6=\frac{m(m-1)(36m^3-54m^2+31m-\frac{13+(-1)^{m-1}5}{2})}{20},$$
$${\rm DT}^{(m)}_7=\frac{7m(m-1)(343m^4-686m^3+539m^2-196m+36)}{720}.$$

\remark In general, ${\rm DT}^{(m)}_n$ has leading term $\frac{n^{n-2}}{n!}m^{n-1}$ considered as a function of $m$. It would be interesting to give a graph-theoretic explanation of this, in the spirit of the graph-theoretic explanation for the leading term of the polynomial counting isomorphism classes of absolutely indecomposable representations of dimension $n$ of the $m$-loop quiver in \cite{HeRV}.\\[2ex]
In \cite{KSCOHA}, a conjecture is formulated which implies the above theorem; we formulate a slight variant of this conjecture.

\begin{conjecture}\label{conjnum}\cite[Section 2.6]{KSCOHA} There exists a product expansion
$$H(q,(-1)^{m-1}t)=\prod_{n\geq 1}\prod_{k\geq 0}\prod_{l\geq 0}(1-q^{k-l}t^n)^{-(-1)^{(m-1)n}c_{n,k}}$$
for nonnegative integers $c_{n,k}$, such that only finitely many $c_{n,k}$ are nonzero for any fixed $n$.
\end{conjecture}

Assuming this conjecture, we have
$$F(q,(-1)^{m-1}t)=\prod_{n\geq 1}\prod_{k\geq 0}\prod_{l=0}^{n-1}(1-q^{k-l}t^n)^{-(-1)^{(m-1)n}c_{n,k}}$$
and thus
$$F((-1)^{m-1}t)=\prod_{n\geq 1}(1-t^n)^{-(-1)^{(m-1)n}n\sum_kc_{n,k}}.$$
Thus, setting ${\rm DT_n}^{(m)}(q)=\sum_{k\geq 0}c_{n,k}q^k$, the conjecture implies that ${\rm DT}_n^{(m)}(q)$ is a polynomial with nonnegative coefficients, such that ${\rm DT}_n^{(m)}(1)={\rm DT}_n^{(m)}$.\\[1ex]
In the following, we will use a simplified notation for product expansions as in the conjecture, using the $\lambda$-ring exponential ${\rm Exp}$, see Section \ref{appendix}. Using Lemma \ref{Exp}, the product of the conjecture can be rewritten as
$${\rm Exp}(\frac{1}{1-q^{-1}}\sum_{n\geq 1}{\rm DT}_n^{(m)}(q)((-1)^{m-1}t)).$$

\section{The Cohomological Hall algebra}\label{coha}

In this section, we review the definition and the main properties of the Cohomological Hall algebra of \cite{KSCOHA} for the $m$-loop quiver. In particular, we formulate the main conjecture of \cite{KSCOHA} on these algebras and relate it to the conjecture of the previous section.\\[2ex]
For a vector space $V$, we denote by $E_V={\rm End}(V)^m$ the space of $m$-tuples of endomorphisms of $V$. The group $G_V={\rm GL}(V)$ acts on $E_V$ by simultaneous conjugation.
For complex vector spaces $V$ and $W$ of dimension $n_1$ and $n_2$, respectively, we consider the subspace $E_{V,W}$ of $E_{V\oplus W}$ of $m$-tuples of endomorphisms $(\varphi_1,\ldots,\varphi_m)$ respecting the subspace $V$ of $V\oplus W$, that is, such that $\varphi_i(V)\subset V$ for all $i=1,\ldots,m$.
We have an obvious projection map $p:E_{V,W}\rightarrow E_V\times E_W$ mapping $(\varphi_1,\ldots,\varphi_m)$ to $((\varphi_1|_V,\ldots,\varphi_m|_V),(\overline{\varphi_1},\ldots,\overline{\varphi_m})),$ where $\overline{\varphi_i}$ denotes the endomorphism of $W$ induced by $\varphi_i$.
The action of $G_{V\oplus W}$ on $E_{V\oplus W}$ restricts to an action of the parabolic subgroup $P_{V,W}$ of $G_{V\oplus W}$, consisting of automorphisms respecting the subspace $V$, on $E_{V,W}$. The projection $p$ is equivariant, if the action of $P_{V,W}$ on $E_V\times E_W$ is defined through the Levi subgroup $G_V\times G_W$ of $P_{V,W}$. Moreover, the closed embedding of $E_{V,W}$ into $E_{V\oplus W}$ is $P_{V,W}$-equivariant.\\[1ex]
Using these maps $E_V\times E_W\leftarrow E_{V,W}\rightarrow E_{V\oplus W}$ and their $P_{V,W}$-equivariance, we can define the following map in equivariant cohomology with rational coefficients:

%We have a projection map $\pi:G_{V\oplus W}\times^{P_{V,W}}E_{V,W}\rightarrow E_{V\oplus W}$ induced by the action of $G_{V\oplus W}$ on $E_{V\oplus W}$.

\begin{eqnarray*}H^*_{G_V}(E_V)\otimes H^*_{G_W}(E_W)&\simeq&H^*_{G_V\times G_W}(E_V\times E_W)\\
&\simeq&H^*_{P_{V,W}}(E_{V,W})\\
&\rightarrow&H^{*+2s_1}_{P_{V,W}}(E_{V\oplus W})\\
%&\simeq&H^*_{G_{V\oplus W}}(G_{V\oplus W}\times^{P_{V,W}}E_{V,W})\\
&\rightarrow&H^{*+2s_1+2s_2}_{G_{V\oplus W}}(E_{V\oplus W}),\end{eqnarray*}
where the shifts in cohomological degree are $s_1=\dim E_{V\oplus W}-\dim E_{V,W}=m\dim V\dim W$ and $s_2=-\dim G_{V\oplus W}/P_{V,W}=-\dim V\dim W$ (see \cite[Section 2.2.]{KSCOHA} for the details). Then the following holds:

\begin{theorem}\cite[Theorem 1]{KSCOHA} The above maps induce an associative unital ${\bf Q}$-algebra structure on $\mathcal{H}=\bigoplus_{n\geq 0}H^{*}_{G_{{\bf C}^n}}(E_{{\bf C}^n})$, which is ${\bf N}\times{\bf Z}$-bigraded if $H^{k}_{G_{{\bf C}^n}}(E_{{\bf C}^n})$ is placed in bidegree $(n,(m-1)\binom{n}{2}-k/2)$.
\end{theorem}

The algebra $\mathcal{H}$ is called the Cohomological Hall algebra of the $m$-loop quiver in\cite{KSCOHA}. Note that, although the vector space underlying $\mathcal{H}$ is independent of $m$ (since all spaces $E_V$ are contractible), the algebra structure depends on $m$.\\[1ex]
The above bigrading differs slightly from the one in \cite{KSCOHA}; it is more suited to our purposes of studying the series $H(q,t)$ in relation to the generating series $F(q,t)$ of Poincar\'e polynomials of ${\rm Hilb}_d^{(m)}$.\\[1ex]
We consider the Poincar\'e-Hilbert series of $\mathcal{H}$: $$P_{\mathcal{H}}(q,t)=\sum_{n\geq 0}\sum_{k\in{\bf Z}}\dim_{\bf Q}\mathcal{H}_{n,k}q^kt^n.$$

\begin{lemma} The series $P_\mathcal{H}(q,t)$ equals $H(q,t)$.
\end{lemma}

\proof The homogeneous component of $\mathcal{H}$ with respect to the first component of the bidegree equals $H^*_{G_{{\bf C}^n}}(E_{{\bf C}^n})\simeq H^*_{G_{{\bf C}^n}}({\rm pt})$, which is isomorphic to a polynomial ring in $n$ generators placed in bidegree $(n,(m-1)\binom{n}{2}-i)$ for $i=1,\ldots,n$. Thus, this component has Poincar\'e-Hilbert series
$$\frac{q^{(m-1)\binom{n}{2}}t^n}{(1-q^{-1})\cdot\ldots\cdot(1-q^{-n})}.$$\hb

Using torus fixed point localization, one obtains the following algebraic description of $\mathcal{H}$:

\begin{theorem}\label{cohasymmetric}\cite[Theorem 2]{KSCOHA} The algebra $\mathcal{H}$ is isomorphic to the following shuffle-type algebra structure on $\bigoplus_{n\geq 0}{\bf Q}[x_1,\ldots,x_n]^{S_n}$, the space of symmetric polynomials in all possible numbers of variables: $$(f_1*f_2)(x_1,\ldots,x_{n_1+n_2})=$$
$$\sum f_1(x_{i_1},\ldots,x_{i_{n_1}})f_2(x_{j_1},\ldots,x_{j_{n_2}})(\prod_{k=1}^{n_1}\prod_{l=1}^{n_2}(x_{j_l}-x_{i_k}))^{m-1},$$
the sum ranging over all shuffles $\{i_1<\ldots<i_{n_1}\}\cup\{j_1<\ldots< j_{n_2}\}=\{1,\ldots,n_1+n_2\}$. A homogeneous symmetric function of degree $k$ in $n$ variables is placed in bidegree $(n,(m-1)\binom{n}{2}-k)$.
\end{theorem} 

From this description we see that $\mathcal{H}$ is commutative in case $m$ is odd, and supercommutative in case $m$ is even.

\begin{conjecture}\label{conjalg}\cite[Conjecture 1]{KSCOHA} The bigraded algebra $\mathcal{H}$ is isomorphic to ${\rm Sym}(C\otimes{\bf Q}[z])$, the (graded) symmetric algebra over a bigraded supervectorspace. For any fixed $n\geq 1$, only finitely many homogeneous components $C_{n,k}$ are nonvanishing and $k\geq 0$ in this case, and $z$ is a homogeneous element of bidegree $(0,-1)$.
\end{conjecture}

Thus conjecture immediately implies Conjecture \ref{conjnum} for $c_{n,k}=\dim_{\bf Q}C_{n,k}$, since the Poincar\'e-Hilbert series of a symmetric algebra has a natural product expansion, namely $P_{{\rm Sym}(V)}={\rm Exp}(P_V)$.

%Fix $m\geq 1$. This is no loss of generality for studying the series $H^{(m)}(q,t)$ due to the following duality result:

%\begin{lemma} We have $H^{(1-m)}(q,t)=H^{(m)}(q^{-1},-qt)$.
%\end{lemma}

\section{The degenerate Cohomological Hall Algebra}\label{dcoha}

We introduce a degenerate form of the Cohomological Hall algebra $\mathcal{H}$ and show that it is of purely combinatorial nature. We analyze its structure using the combinatorics of partitions in the set $T_n$ introduced in Section \ref{Hilb}.\\[2ex]
Denote by $\Lambda_n$ the set of partitions $\lambda=(0\leq\lambda_1\leq\ldots\leq\lambda_n)$ of length $l(\lambda)=n$, and denote the disjoint union of all $\Lambda_n$ (for $n\geq 0$) by $\Lambda$. For $N\in{\bf N}$, define $S^N\lambda=(\lambda_1+N,\ldots,\lambda_n+N)$. Define the union $\mu\cup\nu$ of partitions $\mu,\nu\in\Lambda$ as the partition with parts $\mu_1,\ldots,\mu_{l(\mu)},\nu_1,\ldots,\nu_{l(\nu)}$, resorted in ascending order.\\[1ex]
Generalizing the definition in Section \ref{Hilb}, the weight of a partition $\lambda$ is defined as ${\rm wt}(\lambda)=(m-1)\binom{n}{2}-|\lambda|$, where $|\lambda|=\lambda_1+\ldots+\lambda_n$. Note that the generating function of $\Lambda$ by weight and length equals $H(q,t)$.

\begin{definition} Define an algebra structure $*$ on the vector space $A$ with basis elements $\lambda\in\Lambda$ by $$\mu*\nu=\mu \cup S^{(m-1)l(\nu)}\mu$$ for $\mu,\nu\in\Lambda$.
\end{definition}

This multiplication is obviously associative, but non-commutative unless $m=1$. It is easy to verify that this algebra is bigraded by weight and length of partitions, and thus has $H(q,t)$ as its  Poincar\'e series.\\[1ex]
The explicit description of the Cohomological Hall algebra in Theorem \ref{cohasymmetric} allows us to define the following (naive) quantization.

\begin{definition} Define the quantized Cohomological Hall algebra $\mathcal{H}_q$ as the bigraded ${\bf Q}[q]$-module $\bigoplus_{n\geq 0}{\bf Q}[q][x_1,\ldots,x_n]^{S_n}$ with the product $$(f_1*f_2)(x_1,\ldots,x_{n_1+n_2})=$$
$$\sum f_1(x_{i_1},\ldots,x_{i_{n_1}})f_2(x_{j_1},\ldots,x_{j_{n_2}})(\prod_{k=1}^{n_1}\prod_{l=1}^{n_2}(x_{j_l}-qx_{i_k}))^{m-1}.$$

\end{definition}

\remark It would be interesting to realize this algebra geometrically, as the convolution algebra in some appropriate cohomology theory on the $G_V$-spaces $E_V$ of the previous section.\\[2ex]
We can specialize the algebra $\mathcal{H}_q$ to any $q\in{\bf Q}$, in particular to $q=0$, yielding an algebra $\mathcal{H}_0$.

\begin{proposition} We have an isomorphism of bigraded algebras $A\simeq\mathcal{H}_0$ by mapping a partition $\lambda$ to the symmetric polynomial $$P_\lambda(x_1,\ldots,x_n)=\sum_{\sigma\in S_n}x_{\sigma(1)}^{\lambda_1}\cdot\ldots\cdot x_{\sigma(n)}^{\lambda_n}.$$
\end{proposition}

\proof The polynomial $P_\lambda$ is a suitable multiple of the monomial symmetric polynomial $m_\lambda(x_1,\ldots,x_n)$. The multiplication in $\mathcal{H}_0$ reduces to $$(f_1*f_2)(x_1,\ldots,x_{n_1+n_2})=$$
$$\sum f_1(x_{i_1},\ldots,x_{i_{n_1}})f_2(x_{j_1},\ldots,x_{j_{n_2}})(\prod_{l=1}^{n_2}x_{j_l})^{(m-1)n_1}.$$
Identification of shuffles with cosets $S_{n_1+n_2}/(S_{n_1}\times S_{n_2})$ immediately shows that $P_\lambda*P_\mu=P_{\lambda*\mu}$.\hb

Recall from section \ref{Hilb} the subset $T_n\subset\Lambda_n$ of partitions $\lambda\in\Lambda_n$ such that $\lambda_i\leq(m-1)(i-1)$ for all $i=1,\ldots,n$, and define $T$ as the disjoint union of all $T_n$.

\begin{lemma} The subspace $B$ of $A$ generated by the basis elements indexed by $T$ is stable under the multiplication $*$, thus $B$ is a subalgebra of $A$. In computing a product $\lambda*\mu$ for $\lambda,\mu\in T$, it suffices to append $S^{(m-1)l(\lambda)}\mu$ to $\lambda$ (without resorting parts).
\end{lemma}

\proof Using the definition of $T$ and of $*$, this is immediately verified.\hb

Denote by $S$ the linear operator on $A$ induced by the operation $S$ on partitions.

\begin{lemma} Multiplication induces an isomorphism of bigraded vector spaces $B\otimes SA\simeq A$.
\end{lemma}

\proof On the level of partitions, this reduces to the statement that multiplication induces a bijection between $\bigcup_{k+l=n}T_k\times (\Lambda_l+1)$ and $\Lambda_n$ preserving weights. Suppose $\lambda$ is given. If $\lambda\in T_n$, we map $\lambda$ to $(\lambda,())\in T_n\times\Lambda_0$. Otherwise, let $i$ be maximal such that $\lambda_i\leq(m-1)(i-1)$ (thus $i<n$). We define $\mu=(\lambda_1,\ldots,\lambda_i)$. We have $\lambda_j>(m-1)(j-1)$ for all $j>i$, thus $(\lambda_{i+1},\ldots,\lambda_n)=S^{(m-1)i+1}\nu$ for the partition $\nu$ of length $n-i$ with parts $\nu_k=\lambda_{i+k}-(m-1)i-1\geq 0$. Then $\lambda$ is mapped to $(\mu,\nu)$. By a simple calculation, compatibility of this bijection with the weight is verified.\hb

We iterate this lemma and get:

\begin{corollary} Multiplication induces an isomorphism $$\bigotimes_{i\geq 0} S^iB=B\otimes SB\otimes S^2B\otimes\ldots\simeq A.$$
\end{corollary}

\proof Iteration of the previous lemma shows that any $\lambda$ admits a finite decomposition $\lambda=\lambda^1*\ldots*\lambda^s$ such that $\lambda^k\in S^kB$ for degree reasons.\hb

Next, we analyze the structure of the algebra $B$. Denote by $T^0_n\subset T_n$ the subset of all $\lambda\in T_n$ such that $\lambda_i<(m-1)(i-1)$ for $i=2,\ldots,n$, by $T^0$ the disjoint union of all $T^0_n$, and by $B^0$ the subspace of $B$ linearly generated by $T^0$.

\begin{lemma} $B$ is isomorphic to the tensor algebra $T(B^0)$.
\end{lemma}

\proof In a product $\lambda=\lambda^1*\ldots*\lambda^k$ of partitions $\lambda^i\in T^0$, the set of indices $l=2,\ldots,n$ such that $\lambda_l=(m-1)(l-1)$ is precisely the set $\{l(\lambda^1)+1,l(\lambda^1)+l(\lambda^2)+1,\ldots,l(\lambda^1)+\ldots+l(\lambda^{k-1})+1\}$. This observation shows that any $\lambda\in T$ admits a unique such decomposition.\hb

We define a total ordering on $T^0$ by the lexicographic ordering, viewing partitions as words in the alphabet ${\bf N}$. This induces a total ordering, the lexicographic ordering in the alphabet $T^0$, on words in $T^0$. Call a word in the alphabet $T^0$ Lyndon if it is strictly bigger than all its cyclic shifts. Denote by $T^L$ the set of all $\lambda^1*\ldots*\lambda^k$ for $\lambda^1\ldots\lambda^k$ a Lyndon word in $T^0$, thus $T^L$ is the union of all $T^L_n=T^L\cap T_n$, and by $B^L$ the subspace of $B$ generated by $T^L$.

\begin{lemma} Multiplication induces an isomorphism of bigraded vector spaces ${\rm Sym}(B^L)\simeq B$.
\end{lemma}

\proof By the previous lemma, we have $B\simeq T(B^0)$, thus $B\simeq {\rm Sym}(L(B^0))$ as vector spaces by Poincar\'e-Birkhoff-Witt, where $L(B^0)$ is the free Lie algebra in $B^0$ (since the free algebra of a vector space is the enveloping algebra of its free Lie algebra). By general results on free Lie algebras \cite{Reu}, the Lyndon words form a basis of the free Lie algebra, since every word can be written uniquely as a product of Lyndon words, weakly increasing with respect to lexicographic ordering on words.\hb

Combining the above lemmas, we arrive at the following description of the algebra $A$:

\begin{theorem}\label{stra} We have an isomorphism of bigraded vector spaces
$$A\simeq{\rm Sym}(\bigoplus_{i\geq 0}S^iB^L).$$
\end{theorem}

\proof The result follows from the following chain of isomorphisms:
$$A\simeq\bigotimes_{i\geq 0}S^iB\simeq\bigotimes_{i\geq 0}S^i{\rm Sym}(B^L)\simeq\bigotimes_{i\geq 0}{\rm Sym}(S^iB^L)\simeq{\rm Sym}(\bigoplus_{i\geq 0}S^iB^L).$$\hb

Note that this result is not a direct analogue of Conjecture \ref{conjalg} for the algebra $A\simeq\mathcal{H}_0$, since the operator $S$ induces a shift of $(0,-n)$ in bidegree on a homogeneous component $B^L_{(n,k)}$ of $B^L$.\\[2ex]
Comparing Poincar\'e-Hilbert series of both sides in the above description of $A$, we get:

\begin{corollary}\label{corqn} We have the following product expansion:
$$H(q,t)={\rm Exp}(\sum_{n\geq 1}\frac{1}{1-q^{-n}}\sum_{\lambda\in T_n^L}q^{{\rm wt}(\lambda)}t^n)).$$
\end{corollary}

%Defining $\widehat{T}$ as the set of all $\lambda+N$ for all $\lambda\in T$ and $N\in{\bf N}$, and $\widehat{T^0}$, $\widehat{T^L}$ similarly, we can reformulate this calculation as:

%\begin{proposition} We have $$H(q,t)={\rm Exp}(\sum_{\lambda\in\widehat{T^L}}q^{{\rm wt}(\lambda)}t^{l(\lambda)}).$$
%\end{proposition}

For application to (quantized) Donaldson-Thomas invariants, we have to describe $H(q,(-1)^{m-1}t)$, thus it is necessary to derive a signed analogue of the previous proposition. Define $T^{L,+}$ as $T^L$ if $m$ is odd, and as $$T^{L,+}=T^L\cup\{\lambda*\lambda\,|\, \lambda\in T^L,\, l(\lambda)\mbox{ odd}\}$$ if $m$ is even. Define $T^{L,+}_n=T^{L,+}\cap T_n$.

\begin{theorem} We have a product expansion
$$H(q,(-1)^{m-1}t)={\rm Exp}(\sum_{n\geq 1}\frac{1}{1-q^{-n}}\sum_{\lambda\in T^{L,+}_n}q^{{\rm wt}(\lambda)}((-1)^{m-1}t)^n).$$
\end{theorem}

\proof If $m$ is odd, there is nothing to prove, so suppose that $m$ is even. From the identity
$$(1+q^at^b)^{-1}={\rm Exp}(q^{2a}t^{2b}-q^at^b)$$
it follows that $H(q,(-1)^{m-1}t)$ equals
$${\rm Exp}(\sum_{n\geq 1}\frac{1}{1-q^{-n}}\sum_{\lambda\in T^L_n}q^{{\rm wt}(\lambda)}((-1)^{m-1}t)^{n}+\sum_{\substack{{n\geq 1}\\{\mbox{\scriptsize odd}}}}\frac{1}{1-q^{-2n}}\sum_{\lambda\in T^L_n}q^{2{\rm wt}(\lambda)}t^{2n}).$$
Now it remains to note that length and weight double when passing from $\lambda$ to $\lambda*\lambda$, and the claim follows.\hb

Arguing as in Section \ref{dt1}, this implies the following combinatorial interpretation of Donaldson-Thomas invariants.

\begin{corollary} We have ${\rm DT}_n^{(m)}=\frac{1}{n}|T_n^{L,+}|$.
\end{corollary}

Define the polynomials $\overline{Q}_n(q)=\sum_{\lambda\in T_n^L}q^{{\rm wt}(\lambda)}$ and ${Q}_n(q)=\sum_{\lambda\in T_n^{L,+}}q^{{\rm wt}(\lambda)}$; we thus have $Q_n(q)=\overline{Q}_n(q)$ except in case $m$ is even and $n=2\overline{n}$ for odd $\overline{n}$, where $Q_n(q)=\overline{Q}_n(q)+\overline{Q}_{\overline{n}}(q^2)$. We can reformulate the above result as
$$H(q,(-1)^{m-1}t)={\rm Exp}(\sum_{n\geq 1}\frac{1}{1-q^{-n}}Q_n(q)((-1)^{m-1}t)^n).$$

\example To illustrate the classes of partitions $T^0\subset T^L\subset T\subset\Lambda$ defined above, we consider the case $m=2$, $n=4$. The set $T_4$ consists of the $14$ partitions
$$(0000),(0001),(0002),(0003),(0011),(0012),(0013),$$
$$(0022),(0023),(0111),(0112),(0113),(0122),(0123).$$
Five of these belong to $T_4^0$; for the other ones, we have the following decompositions:
$$(0003)=(000)*(0), (0013)=(001)*(0), (0022)=(00)*(00),$$
$$(0023)=(00)*(0)*(0), (0111)=(0)*(000), (0112)=(0)*(001),$$
$$(0113)=(0)*(00)*(0), (0122)=(0)*(0)*(00), (0123)=(0)*(0)*(0)*(0).$$
The lexicographic ordering on $T^0$ gives $(0)<_{\rm lex}(00)<_{\rm lex}(000)<_{\rm lex}(001)$,
thus we have the following eight elements in $T_4^L$:
$$(0000),(0001),(0002),(0003),(0011),(0012),(0013),(0023).$$

\section{Explicit formulas and integrality}\label{integrality}

%Next, we want to prove that $\sum_{\lambda\in\widehat{T_n^{L,+}}}q^{{\rm wt}(\lambda)}$ is of the form $\frac{1}{1-q}P_n(q)$ for a polynomial $P_n(q)\in{\bf Z}[q]$. By the definition of $\widehat{T}$, we have
%$$\sum_{\lambda\in\widehat{T_n^{L,+}}}q^{{\rm wt}(\lambda)}=\frac{1}{1-q^n}\sum_{\lambda\in{T_n^{L,+}}}q^{{\rm wt}(\lambda)},$$
%thus we have to prove that $Q_n(q)=\sum_{\lambda\in{T_n^{L,+}}}q^{{\rm wt}(\lambda)}$ is divisible by $[n]=\frac{1-q^n}{1-q}$. For this, it suffices to prove that $Q_n(\zeta)=0$ for each $n$-th root of unity $\zeta\not=1$.\\[1ex]
Denote by $U_n$ the set of all sequences $(a_1,\ldots,a_n)$ of nonnegative integers which sum up to $(m-1)n$. We consider the natural action of the $n$-element cyclic group $C_n$ on $U_n$ by cyclic shift; call a sequence primitive if it is different from all its proper cyclic shifts. Every non-primitive sequence can be written as the $(n/d)$-fold repetition of a primitive sequence in $U_d$ for $d$ a proper divisor of $n$; we denote the corresponding subset of $U_n$ by $U_n^{d-{\rm prim}}$, and in particular by $U_n^{\rm prim}=U_n^{n-{\rm prim}}$ the subset of primitive sequences. We relate $U_n^{\rm prim}/C_n$, the set of $C_n$-orbits of primitive sequences, to the set $T_n^L$ of the previous section.

\begin{lemma} We have an injective map $\varphi$ from $T_n$ to $U_n$ given by
$$(\lambda_1,\ldots,\lambda_n) \mapsto (\lambda_2-\lambda_1,\lambda_3-\lambda_2,\ldots,\lambda_n-\lambda_{n-1},(m-1)n-\lambda_n).$$
Its converse is given by
$$(a_1,\ldots,a_n)\mapsto(0,a_1,a_1+a_2,\ldots,a_1+\ldots+a_{n-1}).$$
The image of $\varphi$ consists of the sequences $(a_1,\ldots,a_n)$ such that $a_1+\ldots+a_i\leq(m-1)i$ for all $i=1,\ldots,n$.
\end{lemma}

\proof This is immediately verified using the definitions.\hb

Call a sequence $(a_1,\ldots,a_n)$ as above admissible if the condition of the previous lemma is satisfied, that is, if it belongs to the image of $\varphi$.

\begin{lemma}\label{adm} Every cyclic class in $U_n$ contains at least one admissible element.
\end{lemma}

\proof Define an auxilliary sequence $(b_1,\ldots,b_n)$ of integers by $b_i=a_i-(m-1)$; then $\sum_ib_i=0$, and the admissibility condition translates into $\sum_{j=1}^ib_j\leq 0$ for all $i\leq n$. Choose an index $i_0$ such that $b_1+\ldots+b_{i_0}$ is maximal among these partial sums. Then $(a_{i_0+1},\ldots,a_n,a_1,\ldots,a_{i_0})$ is admissible: for $i_0\leq i\leq n$ we have
$$b_{i_0+1}+\ldots+b_i=(b_1+\ldots+b_i)-(b_1+\ldots+b_{i_0})\leq 0.$$
For $i\leq i_0$, we have (since the $b_i$ sum up to $0$):
$$b_{i_0+1}+\ldots+b_n+b_1+\ldots+b_i=(b_1+\ldots+b_i)-(b_1+\ldots+b_{i_0})\leq 0.$$\hb

\begin{proposition}\label{bij} The map $\varphi$ induces a bijection between $T_n^L$ and $U_n^{\rm prim}/C_n$.
\end{proposition}

\proof If $\mu\in T_k$ and $\nu\in T_l$ for $k+l=n$, then $\varphi(\mu*\nu)$ is just the concatenation of the sequences $\varphi(\mu)$ and $\varphi(\nu)$. Thus, $\varphi(\mu*\nu)$ and $\varphi(\nu*\mu)$ are cyclic shifts of each other. Conversely, if a sequence $a\in U_n$ and a proper cyclic shift $a'=(a_{i+1},\ldots,a_n,a_1,\ldots,a_i)$ of $a$ are both admissible, both subsequences $(a_1,\ldots,a_i)$ and $(a_{i+1},\ldots,a_n)$ are admissible. It follows that $a=\varphi(\mu*\nu)$ and $a'=\varphi(\nu*\mu)$ for some $\mu,\nu$.\\[1ex]
We conclude that the restriction of $\varphi$ to $T_n^L$ only maps to primitive classes, and that each such cyclic class is hit precisely once.\hb

Define $U_n^{{\rm prim},+}$ as $U_n^{\rm prim}\cup U_n^{\overline{n}-{\rm prim}}$
%$$U_n^{\rm prim}\cup\{(b_1,\ldots,b_{n/2},b_1,\ldots,b_{n/2})\,:\, (b_1,\ldots,b_{n/2})\in U_{n/2}^{\rm prim}\}$$
if $m$ is even and $n=2\overline{n}\equiv 2\bmod 4$, and as $U_n^{\rm prim}$ otherwise. We have the following variant of the previous proposition:

\begin{corollary} The map $\varphi$ induces a bijection between $T_n^{L,+}$ and $U_n^{{\rm prim},+}/C_n$.
\end{corollary}

Under the above map $\varphi$, the weight ${\rm wt}(\lambda)$ of a partition translates into the function
$${\rm wt}(a_1,\ldots,a_n)=\sum_{i=1}^n(n-i)(m-1-a_i).$$

\begin{lemma} Considered modulo $n$, the function ${\rm wt}$ on $U_n$ is invariant under cyclic shift. In each cyclic class, it assumes its maximum at an admissible element. If $a\in U_n$ is the $\frac{n}{d}$-fold repetition of a sequence $b\in U_d$, then ${\rm wt}(a)=\frac{n}{d}{\rm wt}(b)$.
\end{lemma}

\proof We have 
$${\rm wt}(a_{i+1},\ldots,a_n,a_1,\ldots,a_i)={\rm wt}(a_1,\ldots,a_n)-n((m-1)i-a_1-\ldots,a_i),$$
proving the first two claims. It follows from a direct calculation that the function ${\rm wt}$ is additive with respect to concatenation of sequences as above, proving the third claim.\hb

Defining ${\rm wt}(C)$ for a cyclic class $C\in U_n^{{\rm prim},+}/C_n$ as the maximal weight of sequences in class $C$, we can thus rewrite the polynomial ${Q}_n(q)$ of the previous section as ${Q}_n(q)=\sum_{C\in U_n^{{\rm prim},+}/C_n}q^{{\rm wt}(C)}$. We also derive the identity
$$nQ_n(q)\equiv\sum_{a\in U_n^{{\rm prim},+}}q^{{\rm wt}(a)}\bmod(q^n-1).$$

Define $P_n(q)=\sum_{a\in U_n}q^{{\rm wt}(a)}$. Using again the previous lemma, we have
%We write $U_n$ as the disjoint union over all divisors $d|n$ of the subsets $U_n^{d-{\rm prim}}$ consisting of sequences which are the $\frac{n}{d}$-fold repetition of a primitive sequence $b\in U_d^{\rm prim}$; in particular, we have $U_n^{n-{\rm prim}}=U_n^{\rm prim}$????????????UND FUER D???????????????????. Using the previous lemma, we have
$$P_n(q)=\sum_{d|n}\sum_{a\in U_n^{d-{\rm prim}}}q^{{\rm wt}(a)}=\sum_{d|n}\sum_{b\in U_d^{\rm prim}}q^{\frac{n}{d}{\rm wt}(b)}$$
and thus
$$P_n(q)\equiv\sum_{d|n}^d\overline{Q}_d(q^{\frac{n}{d}})\bmod(q^n-1).$$
By Moebius inversion, this gives
\begin{lemma} We have
$$\overline{Q}_n(q)\equiv\frac{1}{n}\sum_{d|n}\mu(\frac{n}{d})P_d(q^{\frac{n}{d}})\bmod(q^n-1).$$
\end{lemma}

\remark Arguments like the above also appear in the context of the ``cyclic sieving phenomenon'' for Gaussian binomial coefficients, see \cite{RSW}.

\begin{theorem}\label{divis} The polynomial $Q_n(q)$ is divisible by $[n]=1+q+\ldots+q^{n-1}$, and the quotient $\frac{1}{[n]}Q_n(q)$ is a polynomial in ${\bf Z}[q]$.
\end{theorem}

\proof First note that $P_n(q)$ equals the $t^{(m-1)n}$-term in
$$\sum_{a_1,\ldots,a_n\geq 0}q^{\sum_i(n-i)(m-1-a_i)}t^{\sum_ia_i}=\frac{q^{(m-1)\binom{n}{2}}}{\prod_{i=0}^{n-1}(1-q^{-i}t)}.$$
%$$=\frac{q^{(m-1)\binom{n}{2}}}{\sum_{0\leq k\leq n}(-1)^kq^{-\binom{k}{2}}\left[{n\atop k}\right]_{q^{-1}}t^k}.$$

Let $\zeta_n$ be a primitive $n$-th root of unity. Specializing $q$ at an arbitrary $n$-th root of unity $\zeta_n^s$ for $s=1,\ldots,n$, we see that $P_n(\zeta_n^s)$ equals the $t^{(m-1)n}$-term in
$$\frac{\zeta_n^{(m-1)\binom{n}{2}s}}{\prod_{i=0}^{n-1}(1-\zeta_n^{si}t)}=\frac{\zeta_n^{(m-1)\binom{n}{2}s}}{(\prod_{i=0}^{\frac{n}{g}-1}(1-\zeta_n^{si}t))^g}=\frac{\zeta_n^{(m-1)\binom{n}{2}s}}{(1-t^\frac{n}{g})^g}=$$
$$\zeta_n^{(m-1)\binom{n}{2}s}\sum_{k\geq 0}{\binom{k+g-1}{g-1}}t^{\frac{n}{g}k},$$
where $g={\rm gcd}(s,n)$. The term $\zeta_n^{(m-1)\binom{n}{2}s}$ is easily seen to equal the sign $(-1)^{(m-1)(n-1)s}$, thus
$$P_n(\zeta_n^s)=(-1)^{(m-1)(n-1)s}\binom{m{\rm gcd}(s,n)-1}{{\rm gcd}(s,n)-1}.$$
Substituting this into the Moebius inversion formula of the previous lemma, we arrive at
$$\overline{Q}_n(\zeta_n^s)=\frac{1}{n}\sum_{d|n}\mu(\frac{n}{d})(-1)^{(m-1)(n-1)s}\binom{m{\rm gcd}(s,n)-1}{{\rm gcd}(s,n)-1}.$$
In particular, we have $\overline{Q}_n(1)=\sum_{d|n}\mu(\frac{n}{d})\binom{mn-1}{n-1}$. Applying Lemma \ref{numbertheory}, we see that $Q_n(\zeta_n^s)=\overline{Q}_n(\zeta_n^s)=0$ except in case $m$ even, $n$ even, $s=\overline{n}=\frac{n}{2}$ odd, where $\overline{Q}_n(-1)=-\frac{1}{\overline{n}}\sum_{d|\overline{n}}\mu(\frac{\overline{n}}{d})\binom{md-1}{d-1}$. Using the above formula for $\overline{Q}_n(1)$, in this case we thus get $Q_n(-1)=\overline{Q}_n(-1)+\overline{Q}_{\overline{n}}(1)=0$.\\[1ex]
We have proved that $Q_n(\zeta_n^s)=0$ for $s=1,\ldots,n-1$, thus $Q_n(q)\in{\bf Z}[q]$ is divisible in ${\bf Z}[q]$ by all nontrivial cyclotomic polynomials $\Phi_d(q)$ for $1\not=d|n$, and thus $Q_n(q)$ is divisible in ${\bf Z}[q]$ by their product, which equals the polynomial $[n]$.\hb

We thus arrive at the following explicit formulas:

\begin{theorem}\label{mainresult} The following holds for all $m\geq 1$ and all $m\geq 1$:
\begin{enumerate}
\item The quantized Donaldson-Thomas invariant ${\rm DT}_n^{(m)}(q)$ is given by
$${\rm DT}_n^{(m)}(q)=q^{1-n}\frac{1}{[n]}\sum_{C\in U_n^{{\rm prim},+}}q^{{\rm wt}(C)}$$
and is a polynomial with integer coefficients.
\item For any $i=0,\ldots,n-1$, the (unquantized) Donaldson-Thomas invariant ${\rm DT}_n^{(m)}$ equals the number of classes $C\in U_n^{{\rm prim},+}$ of weight ${\rm wt}(C)\equiv i\bmod n$.
\end{enumerate}
\end{theorem}

\proof Using Corollary \ref{corqn} and the definition of ${\rm DT}_n^{(m)}(q)$ of Section \ref{dt1}, the first part follows from Theorem \ref{divis}. The second part follows by comparing coefficients of the polynomials $Q_n(q)$ and ${\rm DT}_n^{(m)}(q)$.\hb

\remark There seems to be no natural weight function $s$ on classes $C\in U_n^{{\rm prim},+}$ of weight ${\rm wt}(C)\equiv i\bmod n$ such that $\sum_Cq^s(C)={\rm DT}_n^{(m)}(q)$.

\section{Relation to Higgs moduli}\label{higgs}

For $n\in{\bf N}$ and $d\in{\bf Z}$, define $H_{n,d}$ as the set of all sequences $(l_1,\ldots,l_n)\in{\bf Z}^n$ with the following properties:
\begin{enumerate}
\item $l_{k+1}-l_k+(m-1)\geq 0$ for all $k=1,\ldots,n-1$,
\item $\sum_{i=1}^nl_i=d$,
\item $\frac{\sum_{i=1}^kd_i}{k}\geq \frac{d}{n}$ for all $k<n$.
\end{enumerate}

These sequences arise as certain fixed points (so-called type $(1,\ldots,1)$-fixed points) in the moduli space of ${\rm SL}_n$-Higgs bundles, for the action of ${\bf C}^*$ scaling the Higgs field; see \cite[Proposition 10.1]{HT}. A relation to Conjecture \ref{conjnum} is hinted at in \cite[Remark 4.4.6]{HRV}.\\[2ex]
\remark Shifting every entry of such a sequence by $1$ defines a bijection $H_{n,d}\simeq H_{n,d+n}$. We also have a duality $H_{n,d}\simeq H_{n,-d}$ by mapping $(l_1,\ldots,l_n)$ to $(-l_n,\ldots,-l_1)$. The elements of $H_{n,0}$ appear in combinatorics as ``score sequences of complete tournaments'' \cite{Lan}.\\[2ex]
%, but no simple form for the generating function of the cardinalities of the $H_{n,0}$ is known. Also the sequence of cardinalities of $\bigcup_{d=0}^{n-1}H_{n,d}$ seems to be a new integer sequence. \\[2ex]
The aim of this section is to prove

\begin{theorem} If $d$ is coprime to $n$, the cardinality of $H_{n,d}$ equals ${\rm DT}_n^{(m)}$.
\end{theorem}

We continue to work with the sets $U_n$, $U_n^{\rm prim}$, $U_n^{\rm prim,+}$ and $U_n^{\rm prim(,+)}/C_n$ of the previous section. We define a map $\Phi:H_{n,d}\rightarrow U_n$ by associating to $l_*=(l_1,\ldots,l_n)$ the sequence $(a_1,\ldots,a_n)$ defined by
$$a_k=l_{k+1}-l_k+(m-1) \mbox{ for }k=1,\ldots,n,$$
where we formally set $l_{n+1}=l_1$; obviously this map is injective. It is also compatible with cyclic shifts, from which it follows easily that the image of $\Phi$ consists only of primitive sequences, and that each cyclic class in $U_n^{prim}$ is hit at most once by the image of $\Phi$ (compare the proof of Proposition \ref{bij} ). In other words, $\Phi$ induces an embedding of $H_{n,d}$ into $U_n^{\rm prim}/C_n$. The weight of $\Phi(l_*)$ is easily computed as $nl_1-d$, thus it is congruent to $-d\bmod n$. We want to prove that, conversely, every primitive cyclic class $a_*$ of weight ${\rm wt}(a_*)\equiv -d\bmod n$ belongs to the image of $\Phi$. We first choose an arbitrary element $a_*$ in such a cyclic class and associate to it the integers
$$l_k=l_1+\sum_{i=1}^{k-1}a_i-(m-1)(k-1)\mbox{ where }l_1=\frac{{\rm wt}(a_*)+d}{n}.$$
This sequence does not necessarily belong to $H_{n,d}$; only the first two defining conditions are fulfilled a priori, together with the condition $l_1-l_n+(m-1)\geq 0$. We define $k_0$ as the maximal index $k\in\{0,\ldots,n\}$ where $l_1+\ldots+l_k-\frac{d}{n}k$ reaches its minimum. Then the cyclic shift $(l_{k_0+1},\ldots,l_n,l_1,\ldots,l_{k_0})$ fulfills the third defining condition by definition of $k_0$ (compare the proof of Lemma \ref{adm}), and the first two conditions are still valid. From the definition, it also follows that this defines an inverse map to $\Phi$.\\[1ex]
We have thus proved that $H_{n,d}$ is in bijection to $(U_n^{\rm prim}/C_n)_{-d}$, the subset of $U_n^{\rm prim}/C_n$ of sequences of weight $\equiv -d\bmod n$. For parity reasons, sequences of such a weight cannot be twice a shorter sequence, thus we can replace $U_n^{\rm prim}/C_n$ by $U_n^{\rm prim,+}/C_n$. By Theorem \ref{mainresult}, the cardinality of the latter equals ${\rm DT}_n^{(m)}$.\hb

\remark As in the case of $(U_n^{\rm prim}/C_n)_{-d}$, there seems to be no natural weight function on sequences $(l_1,\ldots,l_n)$ as above which gives ${\rm DT}_n^{(m)}(q)$.

%The generating function by weight of $U_n^{\rm prim,+}/C_n$ is the polynomial $Q_n(q)$, which is divisible by $[n]$. This means that all the sets $(U_n^{\rm prim,+}/C_n)_d$ of cyclic classes of weight $\equiv -d\bmod n$ have the same cardinality, and thus the cardinality of $(U_n^{\rm prim}/C_n)_{-d}$ equals $Q_n(1)/n=R_n(1)$. The theorem follows.

\section{Appendix: $\lambda$-ring exponential and Moebius inversion}\label{appendix}

Let $R$ be the ring ${\bf Z}[q,q^{-1}][[t]]$ of formal power series in $t$ with coefficients being integral Laurent series in $q$, and denote by $R^+$ the ideal of formal series without constant term. Then $\exp$ and $\log$ define mutually inverse isomorphisms between the additive group of $R^+$ and the multiplicative group $1+R^+$ of formal series with constant term $1$.\\[1ex]
We can define (see \cite{Ge}) a $\lambda$-ring structure on $R$ with Adams operations $\psi_i$ for $i\geq 1$ given by $\psi_i(q)=q^i$, $\psi(t)=t^i$; thus, in particular, we have $\psi_1={\rm id}$ and $\psi_i\psi_j=\psi_{ij}$. We define $\Psi=\sum_{i\geq 1}\frac{1}{i}\psi_i$.

\begin{lemma} The operator $\Psi$ is invertible with inverse $\Psi^{-1}=\sum_{i\geq 1}\frac{\mu(i)}{i}\psi_i$, where $\mu$ denotes the number-theoretic M\"obius function.
\end{lemma}

\proof Composition of $\Psi$ with the operator defined on the right hand side of the claimed formula yields $\sum_{n\geq 1}\sum_{i|n}\mu(i)\frac{\psi_n}{n}$, which equals $\psi_1={\rm id}$ by properties of the Moebius function.\hb

Using this explicit form of the operators $\Psi$ and $\Psi^{-1}$, we can derive the following $q$-Moebius inversion formula for polynomials $f_n(q),g_n(q)$ in $q$ (viewing them as coefficients of formal series):

\begin{lemma} We have $$ng_n(q)=\sum_{d|n}f_d(q^{n/d)})\iff f_n(q)=\sum_{d|n}\mu(\frac{n}{d})dg_d(q^{n/d}).$$
\end{lemma}

We define the $\lambda$-ring exponential ${\rm Exp}:R^+\rightarrow 1+R^+$ by ${\rm Exp}=\exp\circ\Psi$. Its inverse is the $\lambda$-ring logarithm ${\rm Log}=\Psi^{-1}\circ \log$. We have the following explicit formula:

\begin{lemma}\label{Exp} For coefficients $c_{i,k}\in{\bf Z}$ such that, for fixed $i\in{\bf N}$, we have $c_{i,k}\not=0$ for only finitely many $k\in{\bf Z}$, the following formula holds:
$${\rm Exp}(\sum_{i\geq 1}\sum_{k\in{\bf Z}}c_{i,k}q^kt^i)=\prod_{i\geq 1}\prod_{k\in{\bf Z}}(1-q^kt^i)^{-c_{i,k}}.$$
\end{lemma}

\proof It suffices to compute ${\rm Exp}(q^kt^i)$, which is
$$\exp(\sum_{j\geq 1}\frac{1}{j}(q^kt^i)^j)=\exp(-\log(1-q^kt^i))=(1-q^kt^i)^{-1}.$$
The lemma follows.\hb

%Recall the $q$-binomal theorem, which we state in the following form:

%\begin{lemma}REFERENZ We have
%$$\sum_{n\geq 0}\frac{t^n}{(1-q^{-1})\cdot\ldots\cdot(1-q^{-n})})=\prod_{n\geq 0}\frac{1}{1-q^{-n}t}.$$
%\end{lemma}

%Using the above notation, this can be rewritten as
%\begin{corollary} $$\sum_{n\geq 0}\frac{t^n}{(1-q^{-1})\cdot\ldots\cdot(1-q^{-n})})={\rm Exp}(\frac{t}{1-q^{-1}}).$$
%\end{corollary}

In Section \ref{integrality}, we make use of the following Moebius inversion type result.

\begin{lemma}\label{numbertheory} Let $f:{\bf N}\rightarrow {\bf Z}$ be function on non-negative integers. For $n\geq 1$ and a proper divisor $s$ of $n$, the sum
$$\frac{1}{n}\sum_{d|n}\mu(\frac{n}{d})(-1)^{(m-1)(d-1)s}f({\rm gcd}(d,s))$$
equals $0$, except when $m$ is even, $n$ is even, $s=\overline{n}=\frac{n}{2}$ is odd, where it equals $-\frac{1}{\overline{n}}\sum_{d|\overline{n}}\mu(\frac{\overline{n}}{d})f(d)$.
\end{lemma}

\proof Suppose first that $m$ is even, $n$ is even, and $s=\overline{n}$ is odd. Every divisor of $n$ is a divisor $d$ of $\overline{n}$ or twice such a $d$. We can then split the sum in question into
$$\frac{1}{n}(\sum_{d|\overline{n}}\mu(\frac{n}{d})(-1)^{d-1}f(d)+\sum_{d|\overline{n}}\mu(\frac{\overline{n}}{d})(-1)^{2d-1}f(d)).$$
Since all divisors $d$ are odd, we have $\mu(\frac{n}{d})=-\mu(\frac{\overline{d}}{n})$, and the sum simplifies to
$$-\frac{2}{n}\sum_{d|\overline{n}}\mu(\frac{\overline{n}}{d})f(d)=-\frac{1}{\overline{n}}\sum_{d|\overline{n}}\mu(\frac{\overline{n}}{d})f(d),$$
as claimed.\\[1ex]
Now suppose that $s$ is an arbitrary proper divisor of $n$, but $s\not=\frac{n}{2}$ in case $m$ is even and $n\equiv 2\bmod 4$. We write
$$\frac{1}{n}\sum_{d|n}\mu(\frac{n}{d})(-1)^{(m-1)(d-1)s}f({\rm gcd}(d,s))=$$
$$\frac{1}{n}\sum_{g|s}\sum_{\substack{{d|\frac{n}{g}}\\ {{\rm gcd}(d,\frac{s}{g})=1}}}\mu(\frac{n}{gd})(-1)^{(m-1)(gd-1)s}f(g).$$
We can uniquely decompose $\frac{n}{g}$ as $n_1n_2$, where $n_1$ collects all prime factors of $\frac{n}{g}$ dividing $\frac{s}{g}$; we then have ${\rm gcd}(n_1,n_2)=1$, and the divisors $d$ of $\frac{n}{g}$ such that ${\rm gcd}(d,\frac{s}{g})=1$ are precisely the divisors of $n_2$. Thus, we can rewrite the above sum as
$$\frac{1}{n}\sum_{g|s}\sum_{d|n_2}\mu(n_1)\mu(\frac{n_2}{d})(-1)^{(m-1)(gd-1)s}f(g)=$$
$$\frac{1}{n}(-1)^{(m-1)s}\sum_{g|s}\mu(n_1)f(g)\sum_{d|n_2}\mu(\frac{n_2}{d})(-1)^{(m-1)gds}.$$
By Moebius inversion, the inner sum, temporarily called $\rho(g)$, equals zero except in the case $n_2=1$, or $n_2=2$ and $(m-1)gs$ is even.\\[1ex]
Now suppose that in the above sum, the summand corresponding to a divisor $g$ of $s$ is non-zero, that is, both $\mu(n_1)$ and $\rho(g)$ are non-zero. First consider the case $n_2=1$, thus $n_1=\frac{n}{g}$ is squarefree. Since $s\not=n$, there exists a prime $p$ dividing $\frac{n}{s}$, and thus also $\frac{n}{g}$. Since $n_2=1$, the prime $p$ also divides $\frac{s}{g}$, thus $p^2$ divides $\frac{n}{g}$, a contradiction. Now consider the case $n_2=2$ and $(m-1)gs$ even, thus $n_1=\frac{n}{2g}$ is squarefree. Again, a prime $p$ dividing $\frac{n}{2g}$ also divides $\frac{s}{g}$. If $p$ is odd, the argument of the first case again yields a contradiction. So suppose that $2$ is the only prime divisor of $\frac{n}{g}$, that $n=2^kn'$ for odd $n'$, and $g=2^ln'$. Then $s=2^{l'}n'$ for some $l'\leq l$, and $\frac{n}{2g}=2^{k-l}$. Since $\frac{n}{2g}$ is squarefree, we have $k=l$ or $k=l+1$. If $k=l$, then $s\geq g=n$, a contradiction. If $k=l+1$, then $g=\frac{n}{2}$, both $s$ and $g$ are odd, and thus $m$ is even. But then $n\equiv 2\bmod 4$, and by assumption $s\not=\frac{n}{2}$, a contradiction.\\[1ex]
Thus we see that no summand above can be non-zero, proving the claim.\hb

\end{document}